\documentclass[12pt]{amsart}
\usepackage{amssymb,latexsym}
\usepackage{enumerate}
\usepackage{hyperref}
\usepackage{fullpage}
\usepackage{tikz}
\usepackage{bbm}

\usepackage{fixme}
\fxsetup{
    status=draft,
    author=,
    layout=margin,
    theme=color
}

\makeatletter \@namedef{subjclassname@2010}{%
  \textup{2010} Mathematics Subject Classification}
\makeatother

\newcounter{thm} \numberwithin{thm}{section}
\newtheorem{Theorem}[thm]{Theorem}

\newtheorem{Lemma}[thm]{Lemma}
\newtheorem{Corollary}[thm]{Corollary}

	\newcommand{\RR}[0]{\mathbb R}

\newcommand{\cB}[0]{\mathcal B}

\newcommand{\cL}[0]{\mathcal L}	
	
\newcommand{\lr}[1]{\left(#1\right)}

\author[B. Hanson]{Brandon Hanson} \address{Brandon Hanson \\University of Maine\\
Orono, ME}
\email{brandon.w.hanson@gmail.com}

\author[O. Roche-Newton]{Oliver Roche-Newton} \address{Oliver Roche-Newton, Johann Radon Institute for Computational and Applied Mathematics\\
Linz, Austria}
\email{o.rochenewton@gmail.com}

\author[S. Senger]{Steven Senger} \address{Steven Senger \\ Missouri State University\\
Springfield, MO}
\email{stevensenger@missouristate.edu}

\date{}

\begin{document}

\baselineskip=17pt

\title{Convexity, Superquadratic Growth, and Dot Products}

\date{}

\begin{abstract}
Let $P \subset \mathbb R^2$ be a point set with cardinality $N$. We give an improved bound for the number of dot products determined by $P$, proving that,
\[
|\{ p \cdot q :p,q \in P \}| \gg N^{2/3+c}.
\]
A crucial ingredient in the proof of this bound is a new superquadratic expander involving products and shifts. We prove that, for any finite set $X \subset \mathbb R$, there exist $z,z' \in X$ such that
\[
\left|\frac{(zX+1)^{(2)}(z'X+1)^{(2)}}{(zX+1)^{(2)}(z'X+1)}\right| \gtrsim |X|^{5/2}.
\]
This is derived from a more general result concerning growth of sets defined via convexity and sum sets, and which can be used to prove several other expanders with better than quadratic growth. The proof develops arguments from \cite{HRNR}, and uses predominantly elementary methods.
\end{abstract}

\maketitle

\section{Introduction}

\subsection{Dot products}

For a finite set $P$ of points in $\mathbb R^2$, let $\Lambda(P)$ denote the set of all dot products determined by $P$, that is
\[
\Lambda(P):= \{ p \cdot q : p,q \in P\},
\]
where the dot product of $p=(p_1,p_2)$ and $q=(q_1,q_2)$ is $p \cdot q= p_1q_1+p_2q_2$. In the spirit of the Erd\H{o}s distinct distance problem, one expects that the set $\Lambda(P)$ should be large for any $P$. The bound
\begin{equation} \label{STbasic}
|\Lambda(P)| \gg |P|^{2/3}
\end{equation}
follows from an application of the Szemer\'{e}di-Trotter Theorem.\footnote{Here and throughout this paper, the notation
 $X\gg Y$, $Y \ll X,$ $X=\Omega(Y)$, and $Y=O(X)$ are all equivalent and mean that $X\geq cY$ for some absolute constant $c>0$. $X \approx Y$ and $X=\Theta (Y)$ denote that both $X \gg Y$ and $X \ll Y$ hold. $X \gg_a Y$ means that the implied constant is no longer absolute, but depends on $a$. We also use the notation $X \gtrsim Y $ and $Y \lesssim X$ to denote that $X \gg Y/(\log_2 Y)^c$ for some absolute constant $c>0$.} The Erd\H{o}s distance problem was resolved up to logarithmic factors in a remarkable paper of Guth and Katz \cite{GK}, which proved that the bound
\begin{equation} \label{GKbound}
    |D(P)| \gg \frac{ |P|}{\log |P|}
\end{equation}
holds for any finite point set $P \subset \mathbb R^2$, where $D(P)$ is the set of Euclidean distances determined by pairs of points in $P$. It is widely believed that a bound similar to \eqref{GKbound} also holds for $\Lambda(P)$. However, in stark contrast to the distance problem, this question remains wide open, and the apparent connection between these two problems may be somewhat misleading. 

For instance, the closely related unit distance conjecture claims that any fixed distance can occur at most $N^{1+\epsilon}$ times among a point set of size $N$. The bound 
\begin{equation} \label{4/3}
|\{ (p,q) \in P \times P : \|p-q\|=1 \}| \ll N^{4/3}
\end{equation}
was established by Spencer, Szemer\'{e}di and Trotter \cite{SST}. This bound is widely believed to be sub-optimal, but there has been no significant improvement to \eqref{4/3} in 37 years. On the other hand, a construction in \cite{IRNR} shows that there exists a point set $P \subset \mathbb R^2$ such that
\[
|\{ (p,q) \in P \times P : p \cdot q=1 \}| \gg  N^{4/3}.
\]

In this paper, we give the following improvement to \eqref{STbasic}.

\begin{Theorem} \label{thm:main}
For any finite set $P \subset \mathbb R^2$
\begin{equation} \label{dpgoal}
|\Lambda(P)| \gtrsim |P|^{\frac{2}{3}+\frac{1}{3057}}.
\end{equation}
\end{Theorem}

We remark that a better bound, with exponent $\frac{2}{3}+\frac{1}{39}$, for the case of skew symmetric bilinear forms was obtained in \cite{IRNR}. This covers, for instance, the set of wedge/cross products determined by a point set, i.e. $p \wedge q = p_1q_2-p_2q_1$.

\subsection{Superquadratic growth}

An essential ingredient in the proof of Theorem \ref{thm:main} is a new sum-product type result which shows that a certain expander involving products and shifts gives superquadratic growth. We will prove the following theorem.

\begin{Theorem} \label{cor:prodshift2}
 Let $X$ be a finite set of positive reals. Then there exist
 $z,z' \in X$ such that
\begin{equation} \label{5/2}
\left|\frac{(zX+1)^{(2)}(z'X+1)^{(2)}}{(zX+1)^{(2)}(z'X+1)}\right| \gtrsim |X|^{5/2}.
\end{equation}
\end{Theorem}

Here, $A^{(k)}$ denotes the $k$-fold product set $\{a_1 \cdots a_k : a_1,\dots, a_k \in A \}$. A key quantitative feature of Theorem \ref{cor:prodshift2} is that the exponent in the right hand side of \eqref{5/2} is strictly larger than $2$. We refer to the set in Theorem \ref{cor:prodshift2} as a \textit{superquadratic expander}. There are many sets defined by a combination of arithmetic operations that, one may suspect, grow in such a way. However, given such a set, it is often difficult to prove that it really is a superquadratic expander, particularly as methods from incidence geometry typically lose their potency in this range. The existence of other superquadratic expanders has been established in \cite{BRNZ}, \cite{R}, \cite{S}, \cite{RNW}.

Theorem \ref{cor:prodshift2} is derived from a more general result, the forthcoming Theorem \ref{MainExpander}. The statement of Theorem \ref{MainExpander} is a little technical, depending on some further definitions, and we do not state it in full yet. However, Theorem \ref{MainExpander} can also be used to derive several other new sum-product results. For instance, we will prove the following result about superquadratic growth under addition for cubes of shifts.

\begin{Corollary} \label{cor:cubes}
Let $A \subset \mathbb R$ be a finite set. Then there exists $a,a' \in A$ such that
\[
|2(a+A)^3 + 2(a'+A)^3 - 2 (a+A)^3 - (a'+A)^3| \gtrsim |A|^{5/2}.
\]
\end{Corollary}

Further applications of Theorem \ref{MainExpander} will be given in the forthcoming Section \ref{sec:other}.

The proof of Theorem \ref{MainExpander} uses similar elementary methods to an earlier paper of Rudnev and the first two authors \cite{HRNR}, which was itself built up from a simple and beautiful observation concerning sumsets of convex sets in a paper of Ruzsa, Shakan, Solymosi and Szemer\'{e}di \cite{RSSS}. Our first step towards proving Theorem \ref{MainExpander} and its consequences is to give a simple proof of the following quadratic expander.

\begin{Theorem} \label{thm:quad}
Let $A \subset \mathbb R$ be a finite set and let $f$ be a strictly convex or concave function. Then there exists $a,a' \in A$ such that
\[
|f(a+A) + f(a'+A) -  f(a+A) | \gg |A|^{2}.
\]
\end{Theorem}
Theorem \ref{thm:quad} gives the optimal order of growth. This can be seen by considering the case
\[
A=\{1,2,\dots,N\} ,\,\,\,\,\,\,\, f(x)=x^2.
\]
A further example illustrating the optimality of Theorem \ref{thm:quad} is given in Section \ref{sec:expanders}.

\subsection{Sketch of the proof of Theorem \ref{thm:main}}

The proof of Theorem \ref{thm:main} is rather long, and so we sketch it here in an attempt to convey some of the main ideas behind the proof to the reader. In this sketch we are using Theorem \ref{cor:prodshift2} as a black box.

In fact, let us begin the sketch by assuming something a little stronger than Theorem \ref{cor:prodshift2}: suppose that there exists $z \in X$ such that
\[
\left|\frac{(zX+1)^{(4)}}{(zX+1)^{(3)}}\right| \gtrsim |X|^{5/2}.
\]
That is, we make the simplifying assumption that $z=z'$ in Theorem \ref{cor:prodshift2}. After several applications of Pl\"{u}nnecke-Ruzsa type inequalities, it would then follow that, for any $Y \subset \mathbb R$,
\begin{equation} \label{2/3+}
|Y(zX+1)| \gtrsim |Y|^{\frac{7}{8}}|X|^{\frac{1}{4}+\frac{1}{128}}.
\end{equation}
Quantitatively, the additional $1/128$ in the exponent is important here. Later in the sketch, we will apply this bound in a case whereby $|Y|=N^{2/3}$ and $|X|=N^{1/3}$, in which case \eqref{2/3+} yields the non-trivial bound
\begin{equation} \label{calc}
|Y(zX+1)| \gtrsim N^{\frac{2}{3} + \frac{1}{384}}.
\end{equation}

To prove Theorem \ref{thm:main}, we begin with an application of the Szemer\'{e}di-Trotter Theorem which tells us that we obtain an exponent better than $2/3$ unless we are in a very specific situation. Roughly speaking, this allows us to assume that $P$, our given set of $N$ points, is supported on exactly $N^{1/3}$ lines through the origin, each containing $N^{2/3}$ points from $P$. Furthermore, since the problem is rotation invariant, we may assume that one of these lines is the $x$-axis. This in fact implies that the set of $x$-coordinates involved in $P$ also has size around $N^{2/3}$, otherwise we have too many dot products between a fixed point from $P$ on the $x$-axis and the remaining elements of $P$. Pigeonholing, we deduce that there is some vertical line $\ell_0$, with equation $x=x_0$, which supports around $N^{1/3}$ points of $P$. Simplifying slightly, let's assume that all of the $N^{1/3}$ lines through the origin contain a point on $P \cap \ell_0$.

Let $S$ denote the set of slopes of these $N^{1/3}$ lines, each containing a point from $P \cap \ell_0$ and $N^{2/3}$ points from $P$ in total. Now apply \eqref{2/3+} for the set $S$. It follows that there is some slope $s \in S$ such that the bound
\[
|Y(sS+1)| \gtrsim |Y|^{\frac{7}{8}}|S|^{\frac{1}{4}+\frac{1}{128}}
\]
holds for any $Y \subset \mathbb R$. However, if we consider the set of dot products between points of $P \cap \{y=sx\}$ and $P \cap \ell_0$, this is equal to the set $x_0Y(sS+1)$, where $Y$ is the set of the $N^{2/3}$ $x$-coordinates of $P \cap \{y=sx\}$. Recalling the calculation in \eqref{calc}, we obtain an exponent better than $2/3$, as required.

If we remove the simplifying assumption that $z=z'$, things become a bit trickier when it comes to constructing the set that plays the role of $Y$. We need to choose a set of $x$-coordinates that works with respect to two lines $y=sx$ and $y=s'x$ simultaneously, and it is not immediately obvious that these lines share many $x$-coordinates. 

However, an extra piece of information comes to our aid here. We may consider the dot products between $P \cap \{y=sx\} $ and $P \cap \{y=s'x\}$. This turns out to be essentially a product set of the two sets of $x$-coordinates. Each set has size $N^{2/3}$, and their product set has approximately the same size, otherwise we are done. We can use this information, and some fairly basic additive combinatorial trickery, to deduce the existence of a good set $Y$ which works simultaneously for both of the fixed lines.

\subsection{Structure of this paper}
We begin Section \ref{sec:expanders} by proving Theorem \ref{thm:quad}. The main goal of the section is to prove our general result on superquadratic growth, which is Theorem \ref{MainExpander}. We then apply this bound to prove our main result on growth of products with shifts, Theorem \ref{cor:prodshift2}. In Section \ref{sec:dp} we use Theorem \ref{cor:prodshift2}, or rather a slightly refined statement containing some extra information, to prove our main result, Theorem \ref{thm:main}. In Section \ref{sec:other}, we give several other applications of Theorem \ref{MainExpander}. This includes the proof of Corollary \ref{cor:cubes}.

\section{Superquadratic expansion} \label{sec:expanders}

\subsection{An elementary proof of a quadratic expander}

The basic argument underlying this section is contained in the forthcoming proof of Theorem \ref{thm:quad}. In fact we record the following more general result, which will prove useful.

\begin{Theorem} \label{thm:quadagain}
Let \[A=\{b_1<\ldots<b_N\}\] be a finite set of real numbers, and let $f$ be a strictly convex or strictly concave function. Suppose that $ h<h'$ are parameters and 
\begin{equation} \label{spacing}
\min\{b_{i+1}-b_i:1\leq i\leq N-1\}\geq h'-h.
\end{equation}
Then
\[
|f(h+A)+f(h'+A)-f(h+A)| \gg |A|^2.
\]
In particular, there exists $a,a' \in A$ such that
\[
|f(a+A)+f(a'+A)-f(a+A)| \gg |A|^2.
\]
\end{Theorem}
\begin{proof}
Observe that the intervals
\[
(b_i+h, b_i+h') : 1 \leq i \leq N
\]
all have length $h'-h$ and do not overlap because of the spacing condition. Therefore, by convexity of $f$, the intervals
\[
(f(b_i+h), f(b_i+h')) : 1 \leq i \leq N
\]
are increasing (or decreasing, if $f$ is concave) in length as $i$ increases, and they also do not overlap.  This enables us to squeeze the smaller intervals into the bigger ones. In particular, for all $1\leq i<j \leq N$ we have
\[
f(b_j+h)<f(b_j+h)+(f(b_i+h') - f(b_i+h)) \leq f(b_j+h').
\]
This shows the existence of at least $j-1$ elements of $f(h+A)+f(h'+A)-f(h+A)$ in the interval $(f(b_j+h), f(b_j+h'))$. Summing over all intervals, it therefore follows that
\[
|f(h+A)+f(h'+A)-f(h+A)| \geq  \sum_{j=1}^N (j-1) \gg N^2.
\]

For the final statement of the theorem, let $a$ and $a'$ be two elements of $A$ such that $a'-a$ is the minimum positive element of $A-A$. 
\end{proof}

As we mentioned in the introduction, Theorem \ref{thm:quad} is optimal in the case when $A=[N]$\footnote{Here we are adopting the standard notation $[N]=\{1,\ldots,N\}$.} and $f(x)=x^2$. Furthermore, if we take $A=\{\log b : b \in B \}$ and $f(x)=e^x$, Theorem \ref{thm:quad} gives the bound
\[
|bB+b'B-bB| \gg |B|^2.
\]
This bound is also optimal, again in the case $B=[N]$. In fact, a closer look at the proof shows that we obtain the bound
\begin{equation} \label{2var}
|bB + (b'-b)B| \gg |B|^2.
\end{equation}
A very similar bound to \eqref{2var} was obtained by Garaev using a different elementary argument. He made the observation that, if $B$ is a set of positive integers with maximal element $b_{max}$ and with $d_{min}=b'-b$ being the smallest positive element in $B-B$, then
\[
|b_{max}B+(b'-b)B| = |B|^2,
\]
as there are no repetitions in this sum set.

There is some flexibility in Theorem \ref{thm:quad}, hinted at in Theorem \ref{thm:quadagain}, as to the choice of the elements $a$ and $a'$. Indeed, we chose them so that $a'-a$ is a minimal positive element of $A-A$, but we could also have taken another close pair of elements, with minimal modification of the proof. In particular, the same proof works if we replace $a,a'$ with a \textit{second nearest consecutive pair} $b,b' \in A$. That is, we choose $b,b'$ to be consecutive elements of $A$ such that
\[
[c<c' \text{ are consecutive and } c'-c < b' -b ] \Rightarrow (c,c')=(a,a').
\]
We then have the following refined version of the statement of Theorem \ref{thm:quadagain}, which will be useful later.

\begin{Theorem} \label{thm:quadagainagain}
Let $A \subset \mathbb R$ be a finite set and let $f$ be a strictly convex or concave function. Let $a<a'$ be a nearest consecutive pair of elements in $A$ and let $b<b'$ be a second nearest consecutive pair. Then
\begin{equation} \label{nearest}
|f(a+A)+f(a'+A)-f(a+A)| \gg |A|^2.
\end{equation}
and \begin{equation} \label{2ndnearest}
|f(b+A)+f(b'+A)-f(b+A)| \gg |A|^2.
\end{equation}
\end{Theorem}

\begin{proof}
The proof of Theorem \ref{thm:quadagain} immediately implies \eqref{nearest}. To prove \eqref{2ndnearest}, define $A' \subset A$ to be the set consisting of every second element of $A$. That is, if $A=\{b_1<\dots<b_N \}$, $A'$ is the subset
\[
A'=\{b_2<a_4<b_{2M} \},
\]
where $M= \lfloor N/2 \rfloor$. Note that,
\[
\min \{ b_{2(i+1)} - b_{2i} : 1 \leq i \leq \lfloor N/2 \rfloor \} \geq b'-b
\]
and so the spacing condition \eqref{spacing} is satisfied. Therefore
\begin{align*}
|f(b+A)+f(b'+A)-f(b+A)|\geq |f(b+A')+f(b'+A')-f(b+A')| \gg |A|^2,
\end{align*}
as required.
\end{proof}

In this paper, we are particularly interested in the expander consisting of products of shifts of products, and so we record the following corollary.

\begin{Corollary} \label{cor:prodshiftquad}
Let $X$ be a set of positive reals. Suppose that $z_1,z_1', z_2,z_2' \in X$ with $z_1 < z_1'$ and $z_2<z_2'$, and with the additional properties 
\begin{itemize}
    \item $\log z_1 < \log z_1'$ is a nearest consecutive pair of elements in $\log X$, and
    \item $\log z_2 < \log z_2'$ is a second nearest consecutive pair of elements in $\log X$.
\end{itemize}
Then
\[
\left | \frac{(z_1X+1)(z_1'X+1)}{(z_1X+1)} \right | \gg |X|^2
\]
and
\[
\left | \frac{(z_2X+1)(z_2'X+1)}{(z_2X+1)} \right | \gg |X|^2.
\]

\end{Corollary}

\begin{proof}
Apply Theorem \ref{thm:quadagainagain} with $A= \log X$ and $f(t)=\log(e^t+1)$.
\end{proof}

\subsection{More variables, more growth}

We would like to do better than Theorem \ref{thm:quadagain} using differences of differences, in the spirit of \cite{HRNR}. We begin to work towards this goal by proving the following lemma.

\begin{Lemma}\label{Expansion2}
Let $X=\{x_1<\ldots<x_N\}$ be a finite set of real numbers and suppose that $F$ is a function which is strictly increasing and strictly convex (or concave) on the interval $[x_1,x_N]$. Then
\[|\lr{(X\times F(X))+(X\times F(X))-(X\times F(X))}\cap[x_1,x_N]\times(F(x_1),F(x_N)]|\gg \lr{\frac{N}{\log N}}^3.\]
\end{Lemma}
\begin{proof}
This argument is taken from \cite{HRNR}. We have
\[(X\times F(X))+(X\times F(X))-(X\times F(X))=(X+X-X)\times(F(X)+F(X)-F(X)).\]
We will estimate the number of elements in the appropriate interval coming from of each factor on the right separately.

Let \[D=\{x_{i+1}-x_i:i\in[N-1]\},\] let \[I_d=\{i\in [N-1]:x_{i+1}-x_i=d\},\]
and let $r(d)=|I_d|$. Since
\[N-1=\sum_{0\leq k\leq \log N}\sum_{\substack{d\in D\\ 2^{k-1}\leq r(d)\leq 2^k}}r(d)\]
there is some integer $L$ and set $D'\subseteq D$, say of size $m$, such that 
\begin{equation}\label{mLBound}
    mL\gg \frac{N}{\log N}, 
\end{equation} and $L/2\leq r(d)\leq L$ for each $d\in D'$. If $d,d'\in D'$ are such that $d\leq d'$ then for each $i\in I_{d'}$ we have
\[x_i+d\in (x_i,x_{i+1}]\cap (X+X-X).\]
It follows that, if $D'=\{d_1<\ldots<d_m\}$, we have \begin{equation}\label{KBound}
    |(X+X-X)\cap (x_1,x_N]|\geq \sum_{j=1}^m j|I_{d_j}|\gg \frac{mN}{\log N}.
\end{equation}

Now fix $d\in D'$ and consider all the intervals \[(F(x_i),F(x_{i+1})]=(F(x_i),F(x_{i}+d)],\ i\in I_d. \] There are at least $L/2$ such intervals, they are visibly disjoint, and each has a distinct length in view of the convexity of $F$. Thus if $F(x_j+d)-F(x_j)<F(x_i+d)-F(x_i)$ for some $i,j\in I_d$, we have
\[F(x_i)+F(x_j+d)-F(x_j)\in (F(x_i,F(x_{i+1})]\cap (F(X)+F(X)-F(X)).\]
and by considering such contributions across all $d\in D'$, we deduce 
\[|(F(X)+F(X)-F(X))\cap (F(x_1),F(x_N)]|\gg \frac{(Lm)^2}{m}.\] Finally, combining the estimates \eqref{mLBound} and \eqref{KBound} concludes the proof.\end{proof}

The following is merely a two-dimensional extension of Theorem \ref{thm:quadagain}, obtained by applying the one-dimensional result in each coordinate independently, and stated in a way that is convenient to our remaining arguments. The lemma extends with no difficulty to arbitrary cartesian products in any dimension.
\begin{Lemma}\label{Expansion1}
Suppose that we have a real numbers $h_1,h_1',h_2,h_2'$ with $h_1<h_1'$ and $h_2 < h_2'$. Suppose that  
\[
A=\{a_1<\ldots<a_N\}
\]
is a finite set of reals
which satisfies the spacing condition \[\min\{a_{i+1}-a_i:i\in[N-1]\}\geq \max \{ h_1'-h_1, h_2'-h_2 \}.
\]
Suppose $f_1$ and $f_2$ are strictly increasing and strictly convex functions. Then, for any $k_1,k_2$ with $N/2\leq k_1,k_2\leq N$, the box  
\[
\cB_{k_1,k_2}=(f_1(a_{k_1}+h_1),f_1(a_{k_1}+h_1')] \times (f_2(a_{k_2}+h_2),f_2(a_{k_2}+h_2')]
\]
contains all points of the form 
\[
(f_1(a_{k_1}+h_1)+f_1(a_{j_1}+h_{1}')-f_1(a_{j_1}+h_1), f_2(a_{k_2}+h_2)+f_2(a_{j_2}+h_{2}')-f_2(a_{j_2}+h_2)),
\] with $1\leq j_1,j_2\leq N/2$. Furthermore, the boxes $\cB_{k_1,k_2}$ are disjoint as $k_1$ and $k_2$ vary.
\end{Lemma}

\begin{proof}
From the spacing condition, we deduce that the boxes
\[
(a_{k_1}+h_1, a_{k_1}+h_1'] \times (a_{k_2}+h_2, a_{k_2}+h_2'] 
\]
are disjoint.
That the boxes $\cB_{k_1,k_2}$ are disjoint follows from this and the fact that each of $f_1$ and $f_2$ is strictly increasing.
To conclude the proof, note that for each $l\in \{1,2\}$ we have $j_l\leq k_l$. The fact that $f_l$ is strictly increasing and strictly convex yields
\[0<f_l(a_{j_l}+h_l')-f_l(a_{j_l}+h_{l})\leq f_l(a_{k_l}+h_l')-f_l(a_{k_l}+h_l).\]
\end{proof}

Lemma \ref{Expansion1} will produce for us a basic amount of expansion, and can be applied rather generally. From there, we will work locally, but in order to find further expansion we will require additional convexity. Such convexity will be found in a pair of discrete derivatives of the function $f$, provided $f$ is itself sufficiently convex. To that end, let $f_1,f_2:[u,\infty)\to\RR$ be continuous functions defined on some interval, which are each strictly increasing and strictly convex. For a non-negative integer $W$, we say the quadruple $(h_1,h_1',h_2,h_2')$ is $W$-good for $(f_1,f_2)$ if $h_1'>h_1>0$, $h_2'>h_2>0$, and if there is a partition 
\[[u,\infty)=I_1\sqcup \cdots\sqcup I_W\]
such that each of the curves 
\[
 \{(f_1(t+h_1')-f_1(t+h_1),f_2(t+h_2')-f_2(t+h_2)) : t\in I_W\}
\] 
is the graph of a strictly convex or strictly concave function. 

This definition just allows us to pass to a large subset of a curve which is convex, should the full curve fail to be so. With it in hand, we are ready to state and prove our main result about superquadratic expanders in its most general form. Again, we remark that the spacing condition serves only to provide additional freedom in our future choice of parameters $h_1,\ h_2,\ h_1',$ and $h_2'$.

\begin{Theorem}\label{MainExpander}
Let $f_1,f_2:[a_1,\infty)\to\RR$ be continuous, strictly increasing and strictly convex functions, and suppose $(h_1,h_1',h_2,h_2')$ is $W$-good for $(f_1,f_2)$, for some $W$. Suppose $A=\{a_1<\ldots<a_N\}$ is a finite set of positive real numbers satisfying the spacing condition
\[\min\{a_{i+1}-a_i:i\in[N-1]\}\geq \max\{h_1'-h_1, h_2'-h_2\}.\] Then
\begin{multline*}|2f_1(A+h_1)-2f_1(A+h_1)+2f_1(A+h_1')-f_1(A+h_1')|\cdot\\\cdot|2f_2(A+h_2)-2f_2(A+h_2)+2f_2(A+h_2')-f_2(A+h_2')| \gg \frac{N^5}{(W\log N)^3}.
\end{multline*}
\end{Theorem}
\begin{proof}
Apply Lemma \ref{Expansion1} to deduce that the disjoint boxes $\cB_{k_1,k_2}$ contain all points of the form
\[(f_1(a_{k_1}+h_1)+f_1(a_{l_1}+h_1')-f_1(a_{l_1}+h_1),f_2(a_{k_2}+h_2)+f_2(a_{l_2}+h_2')-f(a_{l_2}+h_2))\]
whenever $1\leq l_1,l_2\leq N/2\leq k_1,k_2\leq N$.
Observe that for fixed values of $k_1$ and $k_2$, the set of such points can be re-written as the cartesian product
\[(f_1(a_{k_1}+h_1),f_2(a_{k_2}+h_2))+g_1(A')\times g_2(A')\]
where
\[A'=\{a_i:i\in[N/2]\}\] and \[g_j(t)=f_j(t+h_j')-f_j(t+h_j).\]
By $W$-goodness, there is an interval $I$ such that
\[
t\in I\mapsto (f_1(a_{k_1}+h_1),f_2(a_{k_2}+h_2)) + (g_1(t),g_2(t))
\]
 is the graph of a strictly convex (or concave) function, say $F$, and such that $A''=I\cap A'$ has size at least $N/3W$. Define 
 \[
 X=f_1(a_{k_1}+h_1)+g_1(A''),
 \] 
 which is the projection of these points to the $x$-axis, so that 
 \[
 F(X)=f_2(a_{k_2}+h_2)+g_2(A'').
 \] 
 Apply Lemma \ref{Expansion2} to set $X$ and the function $F$, in order to deduce that 
 \[
 |((X\times F(X))+(X\times F(X))-(X\times F(X)))\cap \cB_{k_1,k_2}|\gg \lr{\frac{N}{W\log N}}^3.
 \] 
 We note that
\[
X+X-X\subseteq 2f_1(A+h_1)-2f_1(A+h_1)+2f_1(A+h_1')-f_1(A+h_1')\]
and
\[
F(X)+F(X)-F(X)\subseteq 2f_2(A+h_2)-2f_2(A+h_2)+2f_2(A+h_2')-f_2(A+h_2').
\]
The boxes $\cB_{k_1,k_2}$ are disjoint, so we conclude the proof by adding the contributions from each choice of $k_1$ and $k_2$.
\end{proof}

\subsection{Proof of Theorem \ref{cor:prodshift2}}

We will now use Theorem \ref{MainExpander} to prove Theorem \ref{cor:prodshift2}. We instead prove it in the form of the following statement, which contains some useful extra information.

\begin{Corollary} \label{cor:prodshift2again}
 Let $X$ be a set of positive reals. Suppose that $z_1,z_1', z_2,z_2' \in X$ with $z_1 < z_1'$ and $z_2<z_2'$, and with the additional properties 
\begin{itemize}
    \item $\log z_1 < \log z_1'$ is a nearest consecutive pair of elements in $\log X$, and
    \item $\log z_2 < \log z_2'$ is a second nearest consecutive pair of elements in $\log X$.
\end{itemize} 
Then
\begin{equation} \label{longprod}
 \left|\frac{(z_1X+1)^{(2)}(z_1'X+1)^{(2)}}{(z_1X+1)^{(2)}(z_1'X+1)}\right| \cdot \left|\frac{(z_2X+1)^{(2)}(z_2'X+1)^{(2)}}{(z_2X+1)^{(2)}(z_2'X+1)}\right| 
  \gtrsim |X|^5
\end{equation}
In particular, there exist $z,z' \in X$ such that
\begin{equation}\label{again}
\left|\frac{(zX+1)^{(2)}(z'X+1)^{(2)}}{(zX+1)^{(2)}(z'X+1)}\right| \gtrsim |X|^{5/2}.
\end{equation}
\end{Corollary}
\begin{proof}
Let $A=\log(X)$ and write
\[
A=\{a_1<a_2<\dots a_N\}.
\]
Let $A'$ be the set of all elements of $A$ with even indices, relabelled as $A'=\{b_1,b_2,\dots b_{\lfloor N/2 \rfloor} \}$ and observe that $A'$ satisfies the spacing condition
\[
\min \{b_{i+1}-b_i : i \in [\lfloor N /2 \rfloor -1] \} \geq \max \{\log z_1' - \log z_1, \log z_2' - \log z_2\}= \log (z_2'/z_2).
\]
Define $X'=\exp(A')$. Now take $f_1(t)=f_2(t)=f(t):=\log(e^t+1)$ and observe that 
\begin{align*}
&|2f(A'+\log(z_1))-2f(A'+\log(z_1))+2f(A'+\log(z_1'))-f(A'+\log(z_1'))|
\\&=\left|\frac{(z_1X'+1)^{(2)}(z_1'X'+1)^{(2)}}{(z_1X'+1)^{(2)}(z_1'X'+1)}\right|.
\end{align*}
Similarly,
\begin{align*}
&|2f(A'+\log(z_2))-2f(A'+\log(z_2))+2f(A'+\log(z_2'))-f(A'+\log(z_2'))|
\\&=\left|\frac{(z_2X'+1)^{(2)}(z_2'X'+1)^{(2)}}{(z_2X'+1)^{(2)}(z_2'X'+1)}\right|.
\end{align*} 
With an appropriate value of $W$, we have from Theorem \ref{MainExpander} that
\[
 \left|\frac{(z_1X'+1)^{(2)}(z_1'X'+1)^{(2)}}{(z_1X'+1)^{(2)}(z_1'X'+1)}\right| \cdot \left|\frac{(z_2X'+1)^{(2)}(z_2'X'+1)^{(2)}}{(z_2X'+1)^{(2)}(z_2'X'+1)}\right| \gg  \frac{|X'|^5}{(W\log |X'|)^3}.
\]
We turn to estimating $W$. The curve in question can be parameterized as a shift of
\[\lr{\log(z_1't+1)-\log(z_1t+1),\log(z_2't+1)-\log(z_2t+1)},\,\,\, t>0\]
and so 
\[\frac{dy}{dx}=\frac{z_2'-z_2}{z_1'-z_1}\cdot \frac{z_1z_1't^2 +(z_1+z_1')t+1}{z_2z_2't^2 +(z_2+z_2')t+1}.\]
Therefore
\begin{multline*}
    \frac{d^2y}{dx^2}=\frac{z_2'-z_2}{z_1'-z_1} \cdot\frac{(z_1z_1'(z_2+z_2')-z_2z_2'(z_1+z_1'))t^2+2(z_1z_1'-z_2z_2')t+(z_1+z_1'-z_2-z_2')}{(z_2z_2't^2 +(z_2+z_2')t+1)^2},
\end{multline*}
which is not identically zero for positive values of $z_1,\ z_1',\ z_2$ and $z_2'$ unless $z_1=z_2$ and $z_1'=z_2'$. As long as this is guaranteed, there are at most $4$ points at which this curve can change convexity, and so we can take $W=5$. This completes the proof. 
\end{proof}

We record one more important corollary here, which is that the pair $z,z'$ satisfying \eqref{again} also satisfies the bound
\[
\left | \frac{(zX+1)(z'X+1)}{(zX+1)} \right | \gg |X|^2
\]
by Corollary \ref{cor:prodshiftquad}. We also make a note here of the obvious fact that $z \neq z'$.

\begin{Corollary} \label{cor:combo}
Let $X$ be a set of positive reals. Then there exists $z,z' \in X$, with $z \neq z'$, such that both
\begin{equation} \label{quadcor}
\left | \frac{(zX+1)(z'X+1)}{(zX+1)} \right | \gg |X|^2
\end{equation}
and
\begin{equation}\label{supercor}
\left|\frac{(zX+1)^{(2)}(z'X+1)^{(2)}}{(zX+1)^{(2)}(z'X+1)}\right| \gtrsim |X|^{5/2}.
\end{equation}
hold.

\end{Corollary}

\subsection{The Pl\"{u}nnecke-Ruzsa inequality and applications}

Before we are ready to start the proof of Theorem \ref{thm:main}, we need to manipulate Pl\"{u}nnecke's inequality in order to deduce from Corollary \ref{cor:combo} a sum-product type theorem which is tailored for the dot product problem. We will prove Corollary \ref{cor:SP}, the key feature of which is that we use only product sets (as opposed to including ratio sets) and still obtain better than quadratic growth. We believe it is of independent interest.

We will use the standard Pl\"{u}nnecke-Ruzsa inequality. See Petridis \cite{P} for a modern proof of the statement.

\begin{Theorem} \label{thm:Plun1}
Let $X$ be a finite set in an additive abelian group. Then for any positive integers $k$ and $l$
\[
|kX-lX| \leq \frac{ |X+X|^{k+l}}{|X|^{k+l-1}}
\]
\end{Theorem}

The Ruzsa Triangle Inequality will be used later in the paper. This is the following result.

\begin{Theorem} \label{thm:RTI}
Let $X,Y,Z$ be finite sets in an additive abelian group. Then 
\[
|Y-Z| \leq \frac{|X-Y||X-Z| }{|X|}.
\]
\end{Theorem}
A quick observation is that Theorem \ref{thm:RTI} gives the bound
\begin{equation} \label{RTI2}
|Y-Z|=|(-Y)-(-Z)|\leq \frac{|X-(-Y)||X-(-Z)| }{|X|}=\frac{|X+Y||X+Z| }{|X|}.
\end{equation}

We will also use the following consequence of the Pl\"{u}nnecke-Ruzsa inequality, which was observed by Garaev \cite{G}.

\begin{Theorem} \label{thm:Plun2}
Let $X_1,\dots X_k$ and $Y$ be finite subsets in an additive abelian group. Then
\[
|X_1+\dots  + X_k| \leq \frac{|X_1+Y|\cdots |X_k+Y|}{|Y|^{k-1}}.
\]
\end{Theorem}

We also need a non-trivial bound for the size of $A(A+1)$. The following result of Garaev and Shen \cite{GS} is sufficient for our purposes.

\begin{Theorem} \label{thm:GS} For any finite set $A \subset \mathbb R$ and any non-zero $\lambda \in \mathbb R$,
\[
|A(A+\lambda)| \gg |A|^{5/4}.
\]
\end{Theorem}
We remark that improvements to Theorem \ref{thm:GS} are known, and a recent result of Stevens and Warren \cite{SW} gives the better exponent $49/38$. We use Theorem \ref{thm:GS} in the proof of Theorem \ref{thm:main} in order to simplify the exposition slightly. If we instead apply the result of Stevens and Warren, we obtain a further small improvement, (namely $c=1/2739$ instead of $1/3057$).

We are now ready to derive the following consequence of Corollary \ref{cor:combo}.

\begin{Corollary} \label{cor:SP}
Let $X$ be a set of positive reals. Then there exists $z,z' \in X$ such that 
\[
|(zX+1)^{(4)}(z'X+1)^{(4)}| \gtrsim |X|^{33/16}  .
\]
\end{Corollary}

\begin{proof} Let $z, z' \in X$ be elements given by Corollary \ref{cor:combo} which satisfy the bounds
\begin{equation} \label{quadcoragain}
\left | \frac{(zX+1)(z'X+1)}{(zX+1)} \right | \gg |X|^2
\end{equation}
and
\begin{equation} \label{5/2again}
\left|\frac{(zX+1)^{(2)}(z'X+1)^{(2)}}{(zX+1)^{(2)}(z'X+1)}\right| \gtrsim |X|^{5/2}.
\end{equation}
First we will prove the inequality
\begin{equation} \label{goal1}
|(zX+1)^{(2)}(z'X+1)^{(2)}| \gg |X|^{13/8}   .
\end{equation}
Indeed, let us label
\[
S=(zX+1)(z'X+1).
\]
By \eqref{RTI2} and \eqref{quadcoragain}, we have
\begin{equation} \label{inter}
|X|^2 \ll |S/S| \leq \frac{|SS|^2}{|S|}.
\end{equation}
Furthermore, since $z \neq z'$, it follows from Theorem \ref{thm:GS} that
\[
|S|=\left|\frac{z'}{z}(zX+1)(z'X+1)\right|=\left|\left(z'X+\frac{z'}{z}\right)\left(\left(z'X+\frac{z'}{z}\right)+\left(1-\frac{z'}{z}\right)\right) \right| \gg |X|^{5/4}.
\]
Plugging this into \eqref{inter} and rearranging yields \eqref{goal1}.

Now define
\[
T=(zX+1)^{(2)}(z'X+1)^{(2)}.
\]
By \eqref{RTI2} again and \eqref{5/2again},
\[
|X|^{5/2} \lesssim |T/T| \leq \frac{|TT|^2}{|T|}.
\] 
We can lower bound $|T|$ using \eqref{goal1}. After doing this and rearranging, it follows that
\[
|TT| \gtrsim |X|^{\frac{33}{16}},
\]
as required.

\end{proof}

\section{Dot products} \label{sec:dp}

\subsection{Controlling the number of slopes with the Szemer\'{e}di-Trotter Theorem}

The next lemma uses the Szemer\'{e}di-Trotter Theorem to give a lower bound for the size of $\Lambda(P)$ in terms of the number of lines through the origin. The result is well-known, and can be used to quickly deduce the threshold bound $|\Lambda(P)| \gg |P|^{2/3}$. Moreover, should it be the case that $P$ is incident to at least $|P|^{1/3+c}$ lines through the origin, the exponent $2/3$ can be improved.

\begin{Lemma} \label{lem:dotbasic}
Let $P \subset \mathbb R^2$ be a finite point set. Let $\mathcal L$ denote the set of all lines through the origin incident to $P$. Then there exists $p \in P$ such that
\[
| \{ p \cdot q : q \in P \}| \gg |P|^{1/2}|\mathcal L|^{1/2}.
\]
\end{Lemma}

\begin{proof}
We may assume without loss of generality that each line in $\mathcal L$ has the form $y= \lambda x$, with $\lambda \neq 0$. Otherwise, we can rotate the whole point set so that it avoids the co-ordinate axes, since all of the quantities involved are rotation invariant. This is just done in order to simplify some notation.

Let $\ell_{\lambda}$ denote the line $y=\lambda x$. For each $\ell_{\lambda} \in \mathcal L$, we choose a point $p_{\lambda} \in P$ on the line $\ell_{\lambda}$ arbitrarily. Consider the set 
\[
\{p_{\lambda} \cdot q : q \in P\}
\]
of all dot products determined by $p_{\lambda}$. Note that $p_{\lambda} \cdot q= p_{\lambda} \cdot q'$ if and only if $q$ and $q'$ lie on a line which is perpendicular to $\ell_{\lambda}$. Let $ S_{\lambda}$ denote the set of all lines with slope $-1/\lambda$ which are incident to $P$. It therefore follows that
\begin{equation} \label{equiv}
|S_{\lambda}|=|\{p_{\lambda} \cdot q : q \in P\}|.
\end{equation}
Define $ S$ to be the set of lines
\[
S= \bigcup_{ \ell_{\lambda} \in \mathcal L}  S_{\lambda}.
\]
Note that this union is disjoint, since elements of $ S_{\lambda}$ and $ S_{\lambda'}$ have distinct slopes for $\lambda \neq \lambda'$. By the Szemer\'{e}di-Trotter Theorem
\begin{equation} \label{STapp}
|\mathcal L||P|=\sum_{\ell_{\lambda} \in \mathcal L}I(P,  S_{\lambda})=I(P,  S) \ll |P|^{2/3}|S|^{2/3}+|P|+|S|
\end{equation}
If $|P|$ dominates the right hand side of \eqref{STapp}, then we must have $|\mathcal L| \ll 1$. In this case, there is some line $\ell \in \mathcal L$ such that $|\ell \cap P| \gg |P|$. This set determines $\Omega(|P|)$ dot products, since the set of dot products determined by set of points on a line is in one-to-one correspondence with the product set of the set of magnitudes.

If the first term on the right hand side of \eqref{STapp} is dominant then a rearrangement gives
\begin{equation} \label{linesbound}
|S| \gg |P|^{1/2}|\mathcal L|^{3/2},
\end{equation}
while if $|S|$ dominates the right hand side of \eqref{STapp} then we notice that by definition $|P|\geq|\mathcal L|,$ so we have \[|S| \gg |P||\mathcal L|  \geq |P|^{1/2}|\mathcal L|^{3/2}\] as well, and so \eqref{linesbound} holds in either case. Because $S$ is a disjoint union,
\[
\sum_{\ell_{\lambda} \in \mathcal L}| S_{\lambda}| \gg |P|^{1/2}|\mathcal L|^{3/2} .
\]
It then follows from \eqref{equiv} that there is some $p_{\lambda} \in P$ such that
\[
|\{p_{\lambda} \cdot q : q \in P\}| \gg |P|^{1/2}|\mathcal L|^{1/2}.
\]
\end{proof}

We will take this lemma as a launch point for our theorem. The lemma cannot itself be improved without additional information. Indeed, if $P=[N]\times [N]$ then $\cL$ has size of order $N^2$. Meanwhile, any single point $(m,n)$ makes at most $2N^2$ distinct dot products with the rest of $P$, so this result is optimal. The key will be to use a few different lines from $\cL$, whose slopes will introduce different arithmetic constraints and allow us to invoke an expansion result. This lemma puts us in a regime where this approach is viable.
\subsection{Proof of Theorem \ref{thm:main}}

\begin{proof}
We repeat some notation from the proof of the previous subsection: $\mathcal L$ denotes the set of lines through the origin covering $P$, and $\ell_{\lambda} \in \mathcal L$ is the line with equation $y=\lambda x$. Fix $c=\frac{1}{3057}$. The proof is presented in such a way that $c$ can be viewed as a parameter, and we will calculate the optimal choice for $c$ at the end of the proof.

The idea of the proof is to use the assumption of few dot products in order to deduce arithmetic combinatorial constraints associated to the set of points $P$. We do this by considering dot products arising from various subsets of $P$, each of which is chosen to impose a different constraint. The resulting combination of constraints will combine in such a way as to violate the expansion results proved in the previous section.

\emph{Step 1. There is a large subset, $P'$, of points lying on a small set, $\cL'$, of rich lines through the origin. The $x$-axis is in $\cL'$, and the other lines in $\cL'$ have positive slope.}

Lemma \ref{lem:dotbasic} implies the desired bound \eqref{dpgoal} if
\[
|\mathcal L| \geq |P|^{1/3+2c}.
\]
Therefore, we can henceforth assume that
\begin{equation} \label{assumption}
|\mathcal L| \leq |P|^{1/3+2c}.
\end{equation}


We perform a dyadic decomposition of the line set $\mathcal L$ according to the number of points from $P$ they contain. Note that
\[
\sum_{j=1}^{\lceil \log_2 |P| \rceil } \sum_{\ell \in \mathcal L : 2^{j-1} \leq |\ell \cap P| < 2^j}|\ell \cap P|=\sum_{\ell \in \mathcal L} |\ell \cap P|=|P|.
\]
Therefore, there is a set $\mathcal L' \subset \mathcal L$ and some integer $M$ such that for all $ \ell \in \mathcal L'$, $M \leq |\ell \cap P| \leq 2M$ and 
\begin{equation} \label{prod}
|\mathcal L'| M \gtrsim |P|.
\end{equation}
It then follows from \eqref{assumption} that
\begin{equation} \label{Mbound}
M \gtrsim |P|^{2/3-2c}.
\end{equation}


For the remaining part of the claim, observe first that at least half of the slopes of the lines in $\mathcal L'$ have the same sign. If the most popular sign is negative, we rotate the whole point set by $90$ degrees so that at least half of the lines in $\mathcal L'$ have positive slope. We can then make a further rotation so that the line with smallest positive slope goes to the $x$-axis.
Abusing notation slightly, we redefine the remaining set of lines as $\mathcal L'$. Let $P'$ denote the set of points from $P$ lying on the lines of $\mathcal L'$. We have $|P'| \geq M| \mathcal L'| \gtrsim |P|$. 

\emph{Step 2. There is a large subset, $\cL''$, of lines which intersect the vertical line $x=x_0$ in $P'$.}

Fix a point $p \in P' $ lying on the $x$-axis. The set of dot products
\[
\{ p \cdot q : q \in P'\}
\]
is in bijection with $L_0$, the set of vertical lines incident to $P'$. If we denote by $X$  the set of all $x$-coordinates appearing in $P'$, then $|X|=|L_0|$, and we may therefore assume that $|X| \leq |P|^{2/3+c}$.


For each $x_0 \in X$, define
\[
\mathcal L(x_0) = \{ \ell \in \mathcal L' : \ell \cap \{x=x_0\} \in P'\},
\]
and observe that
\[
\sum_{x_0 \in X} |\mathcal L(x_0)|=|P'| \gtrsim |P|.
\]
It therefore follows by the pigeonhole principle that there exists $x_0 \in X$ such that 
\begin{equation} \label{L0}
|\mathcal L(x_0) | \gtrsim \frac{ |P|}{|X|} \geq |P|^{1/3-c}.
\end{equation}
Let $\mathcal L'':=\mathcal L(x_0)$ be the set satisfying \eqref{L0}. 

\emph{Step 3. A first appeal to expansion --  the set, $S$, of slopes from $\cL''$ contains two elements, $s$ and $s'$, such that \begin{equation} \label{case1}
|(sS+1)^{(4)}(s'S+1)^{(4)}| \gtrsim |S|^{33/16}.
\end{equation}}

This is a direct application of Corollary \ref{cor:SP} (with the role of $X$ played by $S$).

\emph{Step 4. The sets, $X(s)$ and $X(s')$, of $x$-coordinates from $\ell_s\cap P'$ and $\ell_{s
'}\cap P'$ respectively, have multiplicative structure controlled by $\Lambda(P)$. In particular, there is a non-zero $a \in \mathbb R$, for which the intersection \[
Z=X(s) \cap a X(s')
\]
is large:
\begin{equation} \label{Zbound}
|Z| \gtrsim |P|^{2/3 -8c}.
\end{equation}
}


Recall that all lines in $\mathcal L''$ have approximately $M$ points from $P'$ on them, and so
\[
|X(s)|, |X(s')| \approx M.
\]
Now consider the subset of $\Lambda(P)$ consisting of all dot products between points of $P' \cap \{y=sx \}$ and $P' \cap \{y=s'x \}$. This is the set
\[\Lambda_1=
\{(x,sx) \cdot (x',s'x') : x \in X(s), x' \in X(s') \}=X(s)X(s')(1+ss').
\]
Note that $|\Lambda(P)| \geq |\Lambda_1| = |X(s)X(s')|$. Therefore, we may assume that
\begin{equation} \label{smallproduct}
|X(s)X(s')| \leq |P|^{2/3+c}.
\end{equation}
Similarly, by considering only dot products involving points from $P' \cap \{y=sx \}$, we may assume that
\[
|X(s)X(s)| \leq |P|^{2/3+c}.
\]
It therefore follows from the Ruzsa Triangle Inequality, and specifically \eqref{RTI2}, that
\[
|X(s)/X(s')| \leq \frac{|X(s)X(s')||X(s)X(s)|}{|X(s)|} \leq \frac{|P|^{4/3+2c}}{M} .
\]
On the other hand
\[
M^2 \approx |X(s)||X(s')|=\sum_{a \in X(s)/X(s')} |X(s) \cap a X(s')|.
\]
Therefore, by \eqref{Mbound}, there exists $a \in X(s)/X(s')$ such that
\[
|Z|=|X(s) \cap a X(s')| \geq \frac{M^3}{|P|^{4/3+2c}} \gtrsim |P|^{2/3 -8c}.
\]

\emph{Step 5. The multiplicative structure of $Z$ is controlled by $\Lambda(P)$:
\begin{equation} \label{aim}
|(
sS+1)Z|, |(s'S+1)Z| \leq | \Lambda(P)|.
\end{equation}
}

Consider the following subset of $\Lambda(P)$:
\[
\Lambda_2:=\{(x,sx)\cdot (x_0,s_1 x_0) : x \in Z, s_1 \in S \}.
\]
This is indeed a subset of $\Lambda(P)$, since $(x,sx)$ and $(x_0,s_1 x_0)$ are elements of $P$. Furthermore,
\[
|\Lambda_2|=|\{(xx_0(1+ss_1) : x \in Z, s_1 \in S \}|=|Z(1+sS)|.
\]
The first claimed estimate in \eqref{aim} follows immediately. Similarly, the set
\[
\Lambda_3=\{(x,s'x)\cdot (x_0,s_1 x_0) : x \in a^{-1}Z, s_1 \in S \}
\]
is a subset of $\Lambda(P)$ and has cardinality
\[
|\{xx_0(1+s's_1) : x \in a^{-1}Z, s_1 \in S\}|=|Z(1+s'S)|.
\] This shows the second claimed estimate.

\emph{Step 6. A second appeal to expansion -- the sets $(sS+1)Z$ and $(s'S+1)Z$ cannot both be too small.}

Observe first that 
\begin{equation}\label{Sbound}
|S|=|\mathcal L''| \geq |P|^{1/3-c}
\end{equation} 
by \eqref{L0}.

Now apply Theorem \ref{thm:Plun2} in the multiplicative setting with $Y=Z$ and $k=8$, along with the bounds \eqref{case1}, \eqref{Sbound} and \eqref{Zbound}, to obtain
\begin{align*} \label{case11}
|(sS+1)Z|^4|(s'S+1)Z|^4 &\geq |Z|^7|(sS+1)^{(4)}(s'S+1)^{(4)}|\\
&\gtrsim |S|^{33/16}|Z|^7 \\
&\gtrsim |P|^{16/3+1/48-929c/16}.
\end{align*}  

It follows from this and \eqref{aim} that
\[
|\Lambda(P)| \gtrsim |P|^{2/3 + 1/384 - 929c/128}.
\]
We optimise by setting $1/384-929c/128=c$. That is, $c=1/3057$. This concludes the proof.
\end{proof}

\section{Further expander results}
\subsection{Consequences of Theorem \ref{MainExpander}} \label{sec:other}

Theorem \ref{MainExpander} is flexible and provides other types of superquadratic expansion. For instance we can use it to prove Corollary \ref{cor:cubes}. We recall the statement below.

\begin{Corollary} \label{cor:cubesagain}
Let $A \subset \mathbb R$ be a finite set. Then there exists $a,a' \in A$ such that
\[
|2(a+A)^3 + 2(a'+A)^3 - 2 (a+A)^3 - (a'+A)^3| \gtrsim |A|^{5/2}.
\]
\end{Corollary}

\begin{proof}
Let $(h_1,h_1')$ to be a nearest consecutive pair of elements in $A$ and let $(h_2,h_2')$ be a second nearest consecutive pair of elements in $A$. Write $A=\{a_1<a_2<\dots<a_N\}$, and let $B$ be the subset of $A$ consisting of all elements with even index. Relabelling, we have $B=\{b_1<b_2<\dots<b_M\}$, with $M = \lfloor N/2 \rfloor \gg N$.

Now apply Theorem \ref{MainExpander} with this set $B$ and with $f_1(x)=f_2(x)=x^3$. Note that our choices for $B$ and $(h_1,h_1',h_2,h_2')$ imply that the spacing condition of Theorem \ref{MainExpander} is satisfied. We claim that $W=2$ is permissible in the application of the theorem. Assuming the claim, it then follows that 
\begin{multline*}|2(A+h_1)^3-2(A+h_1)^3+2(A+h_1')^3-(A+h_1')^3|\cdot\\\cdot|2(A+h_2)^3-2(A+h_2)^3+2(A+h_2')^3-(A+h_2')^3| \gg \frac{N^5}{(\log N)^3}.
\end{multline*}
Therefore, at least one of the two components of this product is greater $N^{5/2}(\log N)^{-3/2}$, and the proof is complete.

It remains to check that $W=2$ is permissible. The parametric curve we are considering is
\[
((t+h_1')^3- (t+h_1)^3, (t+h_2')^3- (t+h_2)^3).
\]
We calculate that
\[
\frac{d^2y}{dx^2}= \frac{h_2'-h_2}{(h_1'-h_1)^2} \cdot \frac{2(h_1'+h_1-h_2'-h_2)}{3(2t+h_1'+h_1)^3}.
\]
The sign of this expression changes at most once, and therefore the curve can be split into two parts, each of which is strictly convex or concave.


\end{proof}

Theorem \ref{MainExpander} yields several other superquadratic expanders of a similar form, simply by taking $f$ to be a convex function and checking the arguments and corresponding calculations, as we have just done in the proof of Corollary \ref{cor:cubesagain}. For instance, it can be used to prove that there exists $a,a' \in A$ such that
\begin{equation} \label{prodshiftalt}
\left|\frac{(a+A)^{(2)}(a'+A)^{(2)}}{(a+A)^{(2)}(a'+A)} \right | \gtrsim |A|^{5/2}.
\end{equation}
We omit the proof of \eqref{prodshiftalt} as it repeats many of the same calculations in the proof of Corollary \ref{cor:prodshift2again}.

On the other hand, there are some functions where the argument breaks down. For instance, it cannot be used to prove the bound
\[
|2(A+A)^2 + 2(A+A)^2 - 2 (A+A)^2 - (A+A)^2| \gtrsim |A|^{5/2}.
\]
Indeed, this is reassuring, since this bound is not true. If we take $A=\{1,\dots,N\}$ then this is a set of integers contained in the interval $[-12N^2,16N^2]$, and thus has size $O(N^2)$.

On the other hand, we can still use Theorem \ref{MainExpander} to deduce results about growth which concerning sums and differences of $(A+A)^2$. For instance, we can deduce the following result.

\begin{Corollary}
For any $A \subset \mathbb R$, there are elements $a,a'\in A$ such that
\begin{multline} \label{SquaresAndCubes}
|2(A+a)^2-2(A+a)^2+2(A+a')^2-(A+a')^2| \cdot
\\ \cdot |2(A+a)^3-2(A+a)^3+2(A+a')^3-(A+a')^3| \gtrsim |A|^5.
\end{multline}
In particular,
\[
\max \{ |4(A+A)^2-3(A+A)^2|,|4(A+A)^3-3(A+A)^3| \} \gtrsim |A|^{5/2}.
\]
\end{Corollary}
\begin{proof}
We apply Theorem \ref{MainExpander} with $f_1(t)=t^2$ and $f_2(t)=t^3$. This time we take $(h_1,h_1')=(h_2,h_2')$ equal to any pair of nearest neighbours in $A$, which will play the role of $(a,a')$. In this case, the curve in question is parameterized by
\[((t+h_1')^2-(t+h_1)^2,(t+h_1')^3-(t+h_1)^3)\]
which is a parabola, and hence $W=1$ is permissible.
\end{proof}

One can interpret the estimate \eqref{SquaresAndCubes} as saying that, if the set $4(A+A)^2-3(A+A)^2$ attains its minimal possible value $O(N^2)$, then we obtain a cubic order of growth for the set $|4(A+A)^3-3(A+A)^3|$, which is better than that given by Corollary \ref{cor:cubesagain}. The estimate \eqref{SquaresAndCubes} is seen to be best possible in view of the case $X=[N]$.

\subsection{A superquadratic bound for three products and shifts}

We conclude the paper by recording the following consequence of the work in Section \ref{sec:expanders}.

\begin{Corollary} \label{cor:prodshift3}
For any $A \subset \mathbb R$,
\[
|(AA+1)(AA+1)(AA+1)| \gtrsim |A|^{2+\frac{1}{32}}.
\]
\end{Corollary}

In order to prove this, we will need the fact that 
\begin{equation} \label{quadprod}
|(AA+1)(AA+1)| \gtrsim |A|^2
\end{equation}
holds for any $A \subset \mathbb R$. This follows from a simple application of the Szemer\'{e}di-Trotter Theorem. The proof of \eqref{quadprod} is a very small modification of the proof of the main result in \cite{RN}, and also exists implicitly in earlier work of Jones \cite{J}. We give the full statement and proof for completeness.

\begin{Lemma} \label{lem:triples}
For any $A \subset \mathbb R$, there exists $a,a' \in A$ such that
\[
|(aA+1)(a'A+1)| \gg \frac{|A|^2}{\log |A|}.
\]
\end{Lemma}
\begin{proof}
Consider the sum
\begin{equation} \label{Esum}
\sum_{a,a' \in A} E^*(aA+1,a'A+1),
\end{equation}
where $E^*(X,Y)$ denotes the usual multiplicative energy between two sets, that is, the number of solutions to the equation
\[
x_1y_1=x_2y_2 ,\,\,\,\,\,x_i \in X, y_i \in Y.
\]
The quantity in \eqref{Esum} is equal to the number of solutions to the equation
\[
(b_1+a^{-1})(b_2+a'^{-1})=(b_3+a^{-1})(b_4+a'^{-1}), \,\,\,\, a,a',b_1,b_2,b_3,b_4 \in A.
\]
Each solution corresponds to a collinear triple $(-a'^{-1},-a^{-1}), (b_4,b_1), (b_2,b_3)$ of three elements of $(A \cup -A^{-1}) \times (A \cup -A^{-1}) $. Therefore, by the Szemer\'{e}di-Trotter Theorem,
\[
\sum_{a,a' \in A} E^*(aA+1,a'A+1) \ll |A|^4 \log |A|.
\]
By pigeonholing, there exists some $a,a' \in A$ such that $E^*(aA+1,a'A+1) \ll |A|^2\log|A|$. An application of the ubiquitous Cauchy-Schwarz bound
\[
E^*(X,Y) \geq \frac{|X|^2|Y|^2}{|XY|}
\]
completes the proof.

\end{proof}

\begin{proof}[Proof of Corollary \ref{cor:prodshift3}]
By \eqref{again} in the statement of Corollary \ref{cor:prodshift2again}, we have
\[
\left|\frac{(AA+1)(AA+1)(AA+1)(AA+1)}{(AA+1)(AA+1)(AA+1)(AA+1)}\right| \gtrsim |A|^{5/2}.
\]
By Theorem \ref{thm:Plun1} and then Theorem \ref{thm:Plun2}, it follows that
\begin{align*}
    |A|^{5/2} & \lesssim \left|\frac{(AA+1)(AA+1)(AA+1)(AA+1)}{(AA+1)(AA+1)(AA+1)(AA+1)}\right|
    \\ & \leq \frac{|(AA+1)(AA+1)(AA+1)(AA+1)|^4}{|(AA+1)(AA+1)|^3}
    \\& \leq \frac{|(AA+1)(AA+1)(AA+1)|^{16}}{|(AA+1)(AA+1)|^{15}}.
\end{align*}
Rearranging and applying Lemma \ref{lem:triples} completes the proof.

\end{proof}

A very similar argument, instead using the bound \eqref{prodshiftalt} as a starting point, can be used to deduce the bound
\[
|(A+A)(A+A)(A+A)| \gtrsim |A|^{2+\frac{1}{32}}.
\]
This improves the bound on the bound $|(A+A)(A+A)(A+A)| \gtrsim |A|^{2+\frac{1}{392}}$, which was proven in \cite{RNS}. We omit the details to avoid repetition.

\subsection*{Acknowledgements} Oliver Roche-Newton was supported by the Austrian Science Fund FWF Projects P 30405-N32 and P 34180. Brandon Hanson was supported by NSF grant 2001622. We are grateful to Ali Mohammadi, Misha Rudnev and Ilya Shkredov for various helpful discussions. Special thanks go to Audie Warren for patiently listening to the details of the paper during its unstable construction phase.

\end{document}